# ADAPTIVE GOODNESS-OF-FIT TESTS BASED ON SIGNED RANKS[1]

BY ANGELIKA ROHDE

*Weierstraß-Institut Berlin*

Within the nonparametric regression model with unknown regression function $l$ and independent, symmetric errors, a new multiscale signed rank statistic is introduced and a conditional multiple test of the simple hypothesis $l = 0$ against a nonparametric alternative is proposed. This test is distribution-free and exact for finite samples even in the heteroscedastic case. It adapts in a certain sense to the unknown smoothness of the regression function under the alternative, and it is uniformly consistent against alternatives whose sup-norm tends to zero at the fastest possible rate. The test is shown to be asymptotically optimal in two senses: It is rate-optimal adaptive against Hölder classes. Furthermore, its relative asymptotic efficiency with respect to an asymptotically minimax optimal test under sup-norm loss is close to 1 in case of homoscedastic Gaussian errors within a broad range of Hölder classes simultaneously.

**1. Introduction.** Consider the nonparametric regression model with $n$ independent observations

$$Y_i = l(X_i) + \varepsilon_i, \qquad i = 1, \ldots, n,$$

some unknown regression function $l$ on the unit interval and design points $0 \leq X_1 < \cdots < X_n \leq 1$. Throughout this paper, the errors are assumed to be independent and symmetrically distributed around zero, which in particular includes the heteroscedastic case. We postulate Lebesgue continuous error distributions in addition for the sake of simplicity. Within this model, we are interested in identifying subintervals in the design space where $l$ deviates significantly from some hypothetical regression curve $l_o$. For this aim, we develop an exact multiple test of the simple hypothesis "$l = l_o$" against a

Received September 2006; revised April 2007.

[1]Supported in part by the Swiss National Science Foundation.

nonparametric alternative. The method does not require a priori knowledge of the explicit error distributions, and it provides simultaneous confidence statements about deviations of $l$ from $l_o$ with given significance level for arbitrary finite sample size.

For the power investigation of our test, we follow the minimax approach introduced by Ingster (1982, 1993), which permits the set of alternatives to consist of an entire smoothness class, separated from the null hypothesis by some distance $\delta_n$ converging to zero. Typically, the distance to the null hypothesis is quantified by some seminorm $\|\cdot\|$. Then for a given significance level $\alpha$ and some positive number $\delta$ the goal is to find a statistical test $\phi$ whose minimal power

$$\inf_{l \in \mathcal{F}\,:\,\|l-l_o\| \geq \delta} \mathbb{E}_l \phi$$

is as large as possible under the constraint that $\mathbb{E}_{l_o}\phi \leq \alpha$. Approximate solutions for this testing problem are known for various classes $\mathcal{F}$ and seminorms $\|\cdot\|$; see, for instance, Ingster (1987, 1993) for the case of $L_p$-norm and Hölder and Sobolev alternatives, Ermakov (1990) for sharp asymptotic results with respect to the $L_2$-norm and Sobolev alternatives and Lepski (1993) and Lepski and Tsybakov (2000) in case of the supremum norm. It is a general problem that the optimal test $\phi$ may depend on $\mathcal{F}$.

In case of an integral norm $\|\cdot\|$, the problem of adaptive (data-driven) testing a simple or parametric hypothesis is investigated, for example, in Eubank and Hart (1992), Ledwina (1994), Ledwina and Kallenberg (1995), Fan (1996), Fan, Zhang and Zhang (2001), Spokoiny (1996, 1998), Hart (1997) and Horowitz and Spokoiny (2001, 2002). The general procedure is to consider simultaneously a family of test statistics corresponding to different values of smoothing parameters, respectively. As Spokoiny (1996) pointed out, the adaptive approach in case of the $L_2$-norm leads necessarily to suboptimal rates by a factor $\log \log n$. In particular, the tests in Fan (1996) and Spokoiny (1996) are based on the maximum of centered and standardized statistics and (up to this constraint) rate-optimal adaptive against a smooth alternative; see also Fan and Huang (2001). For our purpose, the supremum norm seems to be the most adequate distance. Within the continuous-time Gaussian white noise model, Dümbgen and Spokoiny (2001) have shown that in contrast to the $L_2$-case, adaptive testing with respect to sup-norm loss is actually possible without essential loss of efficiency. They propose a test based on the supremum of suitably standardized kernel estimators of the regression function over different locations and over different bandwidths in order to achieve adaptivity. Unfortunately, their testing procedure depends explicitly on homoscedasticity and Gaussian errors or errors with at least sub-Gaussian tails. If these assumptions are violated, the test may lose its exact or even asymptotic validity. Moreover, its asymptotic power can be arbitrarily small.

In the following section, a new multiscale signed rank statistic is introduced and a conditional test of a one-point hypothesis against a nonparametric alternative is developed. In the third section, its asymptotic power is studied in the setting of homoscedastic errors. A lower bound for minimax testing with respect to sup-norm loss is provided, which is explicitly given in terms of Fisher information. The test turns out to be rate-optimal against arbitrary Hölder classes, provided that the Fisher information of the error distribution is finite. Moreover, a lower bound for its relative asymptotic efficiency with respect to an asymptotically minimax optimal test under sup-norm loss is determined, and the classical efficiency bound $3/\pi$ is recovered even over a broad range of Hölder classes simultaneously. A numerical example illustrating our method is presented in Section 4. Possible extensions are briefly discussed in Section 5. All proofs are deferred to Section 6.

For asymptotic investigations, the design variables are supposed to be deterministic and sufficiently regular in the sense of the assumption

(D) There exists a strictly positive and continuous Lebesgue probability density $h$ on $[0,1]$ of finite total variation such that $X_i = H^{-1}(i/n)$, with $H$ the distribution function of $h$.

By substraction of $l_o$ from the observations, we may assume without loss of generality that $l_o = 0$. Depending on the design density $h$, it is then assumed that under the alternative the regression function $l$ belongs to some smoothness class

$$\mathcal{H}_h(\beta, L) := \{l/\sqrt{h} | l \in \mathcal{H}(\beta, L; [0,1])\},$$

where for any interval $I \subset \mathbb{R}$, $\mathcal{H}(\beta, L; I)$ denotes the class of Hölder functions on $I$ with parameters $\beta, L > 0$. In case $0 < \beta \leq 1$,

$$\mathcal{H}(\beta, L; I) := \{f : I \to \mathbb{R} | |f(x) - f(y)| \leq L|x-y|^\beta \text{ for all } x, y \in I\}.$$

If $k < \beta \leq k+1$ for an integer $k \geq 1$, let $\mathcal{H}(\beta, L; I)$ be the set of functions on $I$ that are $k$ times differentiable and whose $k$th derivative belongs to $\mathcal{H}(\beta - k, L; I)$. We also write $\mathcal{H}(\beta, L)$ for $\mathcal{H}(\beta, L; [0,1])$. In particular, $\mathcal{H}_h(\beta, L)$ coincides with $\mathcal{H}(\beta, L)$ for $h(\cdot) = 1$, corresponding to equidistant design points $X_i = i/n$, $i = 1, \ldots, n$.

**2. The multiscale signed rank statistic.** Inspired by the high asymptotic efficiency of Wilcoxon's signed rank test in simple location shift families [see Hájek and Šidák (1967)], the idea is to define a multiscale testing procedure combining suitably standardized local signed rank statistics. The construction is related to the work of Dümbgen (2002), who used local rank statistics for a test of stochastic monotonicity. In the present context it will turn out that the highest asymptotic efficiency is achieved by weighted local signed rank statistics.

For some kernel function $\psi$ on $[0,1]$ to be specified later and any pair $(s,t)$ with $0 \leq s < t \leq 1$, let $\psi_{st}$ be the shifted and rescaled kernel on the interval $[s,t]$, pointwise given by

$$\psi_{st}(x) := \psi\left(\frac{x-s}{t-s}\right).$$

For notational convenience, we simply write $\psi_{jk}$ for $\psi_{X_j X_k}$, $X_j < X_k$. For any $1 \leq j < k \leq n$ let $R_{jk} := (R_{jk}(i))_{i=j}^k$, with $R_{jk}(i)$ the rank of $|Y_i|$ among the $k-j+1$ numbers $|Y_l|$, $l = j, \ldots, k$. Define the local test statistic

$$T_{jk} := \frac{\sum_{i=j}^k \psi_{jk}(X_i) \operatorname{sign}(Y_i) R_{jk}(i)}{\sqrt{\sum_{i=j}^k \psi_{jk}(X_i)^2 R_{jk}(i)^2}} \tag{1}$$

if the denominator is not equal to zero; and set $T_{jk}$ equal to zero otherwise. The law of $T_{jk}$ depends heavily on the unknown error distributions, but under the null hypothesis, the conditional distribution $\mathcal{L}(T_{jk}|R_{jk})$ does not—even in case of heteroscedastic errors. Hence distribution-freeness may be achieved via conditioning on the ranks. Note that the denominator in (1) is the conditional standard deviation of the numerator given $R_{jk}$ under the null hypothesis.

The question is how to combine these single test statistics in an adequate way. The following theorem acts as a motivation for our approach.

THEOREM 1. *Let the test statistic $T_n$ be defined by*

$$T_n := \max_{1 \leq j < k \leq n} \{|T_{jk}| - \sqrt{2 \log(n/(k-j))}\},$$

*based on a continuous kernel $\psi: [0,1] \to \mathbb{R}$ of bounded total variation with $\int \psi(x)\,dx > 0$. Let assumption* (D) *be satisfied. Then in case of independent identically distributed errors,*

$$\mathcal{L}_0(T_n|R_{1n}) \to_{w,P_0} \mathcal{L}(T_0),$$

*where*

$$T_0 := \sup_{0 \leq s < t \leq 1} \left\{ \frac{|\int_s^t \psi_{st}(x)\sqrt{h(x)}\,dW(x)|}{\|\psi_{st}\sqrt{h}\|_2} - \sqrt{2\log(1/(H(t)-H(s)))} \right\},$$

*with $W$ a Brownian motion on the unit interval.*

Here, $\to_{w,P}$ refers to weak convergence in probability. It follows from results in Dümbgen and Spokoiny (2001) that $T_0$ is finite almost surely. The additive correction in the limiting statistic appears as a suitable calibration for taking the supremum. For it is well known that the maximum of $n$ independent $\mathcal{N}(0,1)$-distributed random variables equals $(2\log n)^{1/2} + o_p(1)$ as $n \to \infty$.

For the testing problem as described in this section, we propose the conditional test

$$\phi_\alpha(Y) := \begin{cases} 0, & \text{if } T_n \leq \kappa_\alpha(R), \\ 1, & \text{if } T_n > \kappa_\alpha(R), \end{cases}$$

where $\kappa_\alpha(R) := \arg\min_{C>0}\{\mathbb{P}(T_n \leq C|R) \geq 1-\alpha\}$ denotes the generalized $(1-\alpha)$-quantile of the conditional distribution $T_n|R$ under the null hypothesis. For explicit applications, we determine $\kappa_\alpha(R)$ via Monte Carlo simulations which are easy to implement. This test is distribution-free and keeps the significance level for arbitrary finite sample size also in the heteroscedastic case. Since the test statistic is discrete-valued, exact level $\alpha$ is attained only for certain values $\alpha \in (0,1)$. In order to achieve arbitrary significance levels exactly, the test can be canonically extended to a randomized procedure.

REMARK (*Simultaneous detection of subregions with significant deviation from zero*). The conditional multiscale test may be viewed as a multiple testing procedure. For a given vector of ranks, the corresponding test statistic $T_n$ exceeds the $(1-\alpha)$-significance level if, and only if, the random family

$$\mathcal{D}_\alpha := \{(X_j, X_k) | 1 \leq j < k \leq n; T_{jk} > \sqrt{2\log(n/(k-j))} + \kappa_\alpha(R)\}$$

is nonempty. Hence one may conclude that with confidence $1-\alpha$, the unknown regression function deviates from zero on *every* interval $(X_j, X_k)$ of $\mathcal{D}_\alpha$.

REMARK (*The choice of the kernel function $\psi$*). If the design density is equal to 1, the limit $T_0$ under the null hypothesis as given in Theorem 1 appears as combination of standardized kernel estimators for the regression function in the standard Gaussian white noise model $dY(t) = l(t)\,dt + n^{-1/2}\,dW(t)$, $0 \leq t \leq 1$. With a certain choice of the kernel $\psi$ depending on the class of alternatives, it coincides there with an asymptotically minimax optimal test statistic with respect to the supremum norm of the testing problem "$l=0$" against Hölder alternatives [Dümbgen and Spokoiny (2001)]. This indicates that in the homoscedastic situation, our conditional test may achieve the highest asymptotic efficiency with the same choice of the kernel function. Here, the construction is as follows: For some Hölder alternative $\mathcal{H}(\beta, L)$, let $\gamma_\beta$ be the solution to the following minimization problem:

(2)    Minimize $\|\gamma\|_2$ over all $\gamma \in \mathcal{H}(\beta, 1; \mathbb{R})$    with $\gamma(0) \geq 1$.

It is known that $\gamma_\beta$ is an even function with compact support, say $[-R, R]$, and $\gamma_\beta(0) = 1 > |\gamma_\beta(x)|$ for $x \neq 0$. To be consistent with the notation introduced above, the optimal kernel $\psi_\beta$ on $[0, 1]$ is then pointwise defined by $\psi_\beta(x) = \gamma_\beta(2Rx - R)$. It is worth noting that the solution $\gamma_\beta$ only depends on the first parameter $\beta$ which shows that the procedure is automatically adaptive with respect to the second parameter $L$. In case $0 < \beta \leq 1$, the solution of (2) is given by $\gamma_\beta(x) = I\{|x| \leq 1\}(1 - |x|^\beta)$. For $\beta > 1$ an explicit solution is known only for $\beta = 2$ [Leonov (1999)]. For details on how this function can be constructed numerically, see Donoho (1994) and Leonov (1999).

**3. Asymptotic power and adaptivity.** In this section, the asymptotic power of our test is investigated in case of independent identically distributed errors. The asymptotic power of the above defined conditional test surely depends on the unknown error distribution as well as the design regularity. The subsequent Theorem 2 provides an extension of Lepski and Tsybakov's (2000) lower bound for the nonparametric regression setting with Gaussian errors to general symmetric error distributions with finite Fisher information. Additionally, the result includes the case of non-equidistant design points.

Let $f$ denote the Lebesgue density of the error distribution. In order to formulate the result on the asymptotic lower bound, let us introduce the following assumptions:

(E1) $f$ is strictly positive and absolutely continuous on $\mathbb{R}$ with finite Fisher information

$$I(f) := \int \left(\frac{f'(x)}{f(x)}\right)^2 f(x)\, dx.$$

The required positivity of the error density $f$ in (E1) just ensures that for any $\theta \in \mathbb{R}$, the shifted distribution $\mathcal{L}_\theta(Y_i) = \mathcal{L}(\varepsilon_i + \theta)$ is absolutely continuous with respect to $\mathcal{L}_0(Y_i) = \mathcal{L}(\varepsilon_i)$. Since we are dealing with noncontiguous alternatives, we are in need of a slightly stronger assumption than differentiability in quadratic mean, which would be equivalent to (E1).

(E2) There exists some positive constant $\delta_0$ such that we have the expansion

$$\int \left\{ \left(\frac{f(z+\theta)}{f(z)}\right)^{1+\delta} - 1 \right\} f(z)\, dz = \frac{1}{2}\delta(1+\delta)\theta^2 I(f)(1 + r(\theta, \delta))$$

with a sequence $r(\theta, \delta) = O(1/\log(1/|\theta|))$ for $|\theta| \to 0$, uniformly in $\delta \in (0, \delta_0]$.

EXAMPLES. (i) (*Normal distribution*). If $f$ denotes the Lebesgue density of the $\mathcal{N}(0,\sigma^2)$-distribution, then $I(f)=\sigma^{-2}$ and

$$\int\left\{\left(\frac{f(z+\theta)}{f(z)}\right)^{1+\delta}-1\right\}f(z)\,dz = \exp\left(\{(1+\delta)^2-(1+\delta)\}\frac{\theta^2}{2\sigma^2}\right)-1$$
$$= \frac{1}{2}\delta(1+\delta)\theta^2 I(f)(1+O(\theta^2))$$

for $\delta$ uniformly bounded from above.

(ii) (*Double exponential distribution*). Let $f$ denote the density of the centered double exponential distribution with parameter $\lambda$, that is, $f(z)=2^{-1}\lambda\exp(-\lambda|z|)$. Simple calculations provide the expansion

$$\int\{(f(z+\theta)/f(z))^{1+\delta}-1\}f(z)\,dz = \tfrac{1}{2}\delta(1+\delta)\theta^2\lambda^2(1+O(\theta)),$$

for $\delta$ uniformly bounded from above, where $\lambda^2=I(f)$.

Via Taylor expansion of $(1+x)^{1+\delta}$ up to the second order and the theorem of dominated convergence, assumption (E2) can be verified for several classical error laws, in particular for the *logistic distribution* which is of exceptional interest in the theory of rank tests. For any $J\subset[0,1]$, let $\|\cdot\|_J$ denote the sup-norm restricted on $J$, that is, $\|l\|_J:=\sup_{x\in J}|l(x)|$.

THEOREM 2. *Let $\rho_n:=((\log n)/n)^{\beta/(2\beta+1)}$ and define*

$$d_* := \left(\frac{2L^{1/\beta}}{(2\beta+1)I(f)\|\gamma_\beta\|_2^2}\right)^{\beta/(2\beta+1)}.$$

*Let the assumptions* (D), (E1) *and* (E2) *be satisfied. Then for arbitrary numbers $\varepsilon_n>0$ with $\lim_{n\to\infty}\varepsilon_n=0$ and $\lim_{n\to\infty}(\log n)^{1/2}\varepsilon_n=\infty$ we obtain*

$$\limsup_{n\to\infty}\inf_{\substack{l\in\mathcal{H}_h(\beta,L):\\ \|l\sqrt{h}\|_J\geq(1-\varepsilon_n)d_*\rho_n}}\mathbb{E}_l\phi_n(Y)\leq\alpha$$

*for any fixed nondegenerate interval $J\subset[0,1]$ and arbitrary tests $\phi_n$ at significance level $\leq\alpha$.*

Even in the knowledge of both smoothness parameters $(\beta,L)$ and the explicit error distribution which is unrealistic for many practical purposes, for any test $\phi_n$ of $\{0\}$ at significance level $\alpha$, there exists an alternative $l$ with $\|l\sqrt{h}\|_J\geq(1-\varepsilon_n)d_*\rho_n$ which will not be detected with probability $1-\alpha-o(1)$ or larger. As expected, the smaller the design density in some location, the more difficult it is to detect there a deviation from zero.

The next theorem is about the asymptotic power of the multiscale signed rank test, based on the kernel being the solution to the minimization problem (2). We restrict our attention to Hölder alternatives with smoothness parameter $\beta \leq 1$. Here the resulting kernel $\psi_\beta$ is pointwise given by $\psi_\beta(x) = (1 - |2x - 1|^\beta)$. For $\beta > 1$, an explicit solution of (2) is known for $\beta = 2$ only; see above. For the sake of simplicity, we consider compact subintervals of $(0, 1)$, which can be avoided by the use of suitable boundary kernels similar to those in Lepski and Tsybakov (2000).

THEOREM 3. *Let $\beta \in (0, 1]$. Let $\phi_n^*$ denote the multiscale signed rank test based on the kernel $\psi_\beta$. Assume that the first derivative of the error density exists and is uniformly bounded and integrable. Denote furthermore $\rho_n := ((\log n)/n)^{\beta/(2\beta+1)}$ and*

$$d^* := \left( \frac{2L^{1/\beta}}{(2\beta + 1)12(\int f(y)^2 \, dy)^2 \|\gamma_\beta\|_2^2} \right)^{\beta/(2\beta+1)}.$$

*Let assumption* (D) *be satisfied and suppose that the modulus of continuity of the design density $h$ is decreasing with at least logarithmic rate, that is, $\sup_{|x-y|\leq \delta} |h(x) - h(y)| = O(1/\log(1/\delta))$ as $\delta \to 0$. Then for arbitrary numbers $\varepsilon_n > 0$ with $\lim_{n\to\infty} \varepsilon_n = 0$ and $\lim_{n\to\infty}(\log n)^{1/2}\varepsilon_n = \infty$ we obtain*

$$\liminf_{n\to\infty} \inf_{\substack{l \in \mathcal{H}_h(\beta, L): \\ \|l\sqrt{h}\|_J \geq (1+\varepsilon_n)d^*\rho_n}} \mathbb{P}_l(\phi_n^* = 1) = 1$$

*for any fixed compact interval $J \subset (0, 1)$.*

The theorem says that if the underlying regression line $l$ multiplied by the square root of the design density deviates from $\{0\}$ by at least $(1 + \varepsilon_n)d^*\rho_n$, then the test rejects the null hypothesis with probability close to 1. Note that the testing procedure does not require knowledge of the design density $h$. Via the choice of the optimal kernel function, the test depends on the smoothness parameter $\beta$, but in contrast to the tests proposed by Lepski and Tsybakov (2000) it remains independent of $L$.

RELATIVE ASYMPTOTIC EFFICIENCY. The ratio $(d_*/d^*)^{(2\beta+1)/\beta}$ may be interpreted as lower bound for the relative asymptotic efficiency in the following sense: Let $(\phi_n)$ be a sequence of arbitrary level-$\alpha$ tests for the simple hypothesis $l = 0$. Let $\delta_n > 0$ such that

$$\liminf_{n\to\infty} \inf_{\substack{l \in \mathcal{H}_h(\beta, L): \\ \|l\sqrt{h}\|_J \geq \delta_n}} \mathbb{E}_l \phi_n = \alpha' > \alpha.$$

Let $m(n)$ be (smallest possible) sample sizes such that
$$\inf_{\substack{l \in \mathcal{H}_h(\beta,L): \\ \|l\sqrt{h}\|_J \geq \delta_n}} \mathbb{E}_l \phi^*_{m(n)} \geq \alpha'.$$
Then under the conditions of Theorems 2 and 3,
$$\liminf_{n \to \infty} \frac{n}{m(n)} \geq (d_*/d^*)^{(2\beta+1)/\beta} = 12 \left( \int f(y)^2 \, dy \right)^2 \Big/ I(f).$$
In case of a Gaussian error density $f = \phi_{0,\sigma^2}$, the former bound equals
$$12\sigma^2 \left( \int \phi_{0,\sigma^2}(y)^2 \, dy \right)^2 = \frac{3}{\pi},$$
which is well known from the classical theory for the Wilcoxon test under the assumption of constant alternatives. The existence of optimal tests for arbitrary error densities $f$ is yet an open problem. In case of homoscedastic Gaussian errors, minimax optimal tests are provided by Dümbgen and Spokoiny (2001). Thus one single test has relative asymptotic efficiency close to 1 with respect to an asymptotically minimax optimal test under supnorm loss for arbitrary Hölder alternatives $\mathcal{H}_h(\beta, L); L > 0$. Sharp asymptotic adaptivity is attained in addition over any range of Hölder classes $\mathcal{H}_h(\beta, L); L_1 \leq L \leq L_2$, for some arbitrary constants $0 < L_1 < L_2 < \infty$. This follows from the fact that the approximations in the proof hold uniformly in $L$ as long as $L$ stays uniformly bounded away from 0 and $\infty$.

Sharp asymptotic adaptivity with respect to both parameters, $\beta$ and $L$, is still an open problem. Nevertheless, under the conditions of Theorems 2 and 3 we obtain the following

THEOREM 4 (Rate-optimality). *Let $\phi_n$ be the conditional multiscale signed rank test at level $\alpha \in (0,1)$, based on some positive continuous kernel $\psi$ of bounded total variation with $\int_0^1 \psi(x) \, d(x) = 1$. Then for arbitrary $\beta > 0$, $L > 0$, there exist constants $c(\beta, L, \psi) \geq d^*(\beta, L)$ such that*
$$\liminf_{n \to \infty} \inf_{\substack{l \in \mathcal{H}_h(\beta,L): \\ \|l\sqrt{h}\|_{[0,1]} \geq c(\beta,L,\psi)\rho_n}} \mathbb{P}_l(\phi_n = 1) = 1.$$

ADAPTIVITY. Without the knowledge of the first parameter $\beta$, the test achieves the optimal rate nevertheless. Note that $\phi_n$ depends neither on $\beta$ nor on $L$. The same considerations concerning the proof as indicated above show that if the range of $(\beta, L)$ is restricted to some compact subset $[\beta_1, \beta_2] \times [L_1, L_2] \subset (0, \infty)^2$, $\phi_n$ is rate-adaptive in the usual setting, that is,
$$\liminf_{n \to \infty} \inf_{(\beta,L) \in [\beta_1,\beta_2] \times [L_1,L_2]} \inf_{\substack{l \in \mathcal{H}_h(\beta,L): \\ \|l\sqrt{h}\|_{[0,1]} \geq c(\beta,L,\psi)\rho_n}} \mathbb{P}_l(\phi_n = 1) = 1.$$

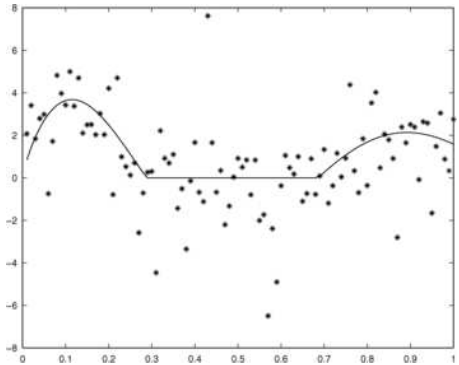

Fig. 1.

REMARK [*Nontrivial power along a sequence of local alternatives* $(l/\sqrt{n})_{n\in\mathbf{N}}$]. In the literature, the power of a goodness-of-fit test is often investigated along a sequence of alternatives $(l/\sqrt{n})_{n\in\mathbf{N}}$. Against such local (but directed) alternatives, the proposed test has nontrivial power as well: If $l$ is continuous with $\|l\|_{\sup} > 0$, then there exists some compact subinterval $J$ of $[0,1]$ with $|l(x)| > \tau > 0$ for all $x \in J$ and some constant $\tau > 0$. The single test statistic $|T_{jk}| - (2\log(n/(k-j)))^{1/2}$ with maximal distance $|X_j - X_k|$ under the constraint $[X_j, X_k] \subset J$ detects a deviation from $\{0\}$ with asymptotic probability arbitrarily close to 1 for sufficiently large $\tau$. Thus, the test is consistent against local alternatives $(a_n l)_{n\in\mathbf{N}}$ whenever $a_n \cdot \sqrt{n} \to \infty$.

**4. Numerical examples.** We illustrate the method with a sample of size $n = 100$ and independent errors drawn from the Student law with three degrees of freedom. The design points are equidistant $X_i = i/n$, and the test statistic is based on the Epanechnikov kernel. Figure 1 shows the regression line with the observations. The estimated quantiles of the conditional test statistic $T_n$ given the vector of ranks of the absolute observation values are based on 999 Monte Carlo simulations. Here we obtained $\kappa_{0.1}(R) = 1.4171$. Figure 2(a) presents the minimal intervals of $\mathcal{D}_{0.1}$, visualized as horizontal line segments and ordered along the $y$-axis in a place-saving manner. Figure 2(b) presents the minimal intervals of rejection at the 0.1-level for an application of the multiscale test [Dümbgen and Spokoiny (2001)], which is based on the idea of homoscedastic Gaussian errors [the standardization by $\sqrt{3} = \mathrm{Var}(\mathrm{Student}_3)^{1/2}$ included]. Based on 999 Monte Carlo simulations as well, we found $\kappa_{0.1} = 1.8187$. The procedure detects a wrong region $[0.56, 0.6]$.

## 5. Extensions.

5.1. *Parametric hypotheses.* Suppose that the null hypothesis $l \in \{l_\theta | \theta \in \Theta\}$ for some parameter space $\Theta \subset \mathbb{R}^d$. If $\widehat{\theta}_n$ denotes a $\sqrt{n}$-consistent estimator of the unknown parameter, the above described procedure is supposed to be applied to the vector of residuals, $(Y_i - l_{\widehat{\theta}_n}(X_i))_{i=1}^n$. In case of equidistant design points and the rectangular kernel, we conjecture that under sufficient regularity conditions on $\widehat{\theta}_n$ and the parametric model, the limit under the null hypothesis of Theorem 1 has the form

$$T_0 := \sup_{0 \leq s < t \leq 1} \left\{ \frac{|W(t) - W(s) + (g(t) - g(s))'Z|}{\sqrt{t-s}} - \sqrt{2\log(1/(t-s))} \right\},$$

with $W$ a Brownian motion on the unit interval, some continuous $\mathbb{R}^d$-valued function $g$ and $Z$ a $d$-variate standard normally distributed random vector. $Z$ comes in via linear expansion of $\widehat{\theta}_n$. The additional estimation of the parameter does not influence the additive correction. However, it destroys the finite sample validity of the conditional test, and a bootstrap procedure may be applied as an approximation.

5.2. *Sobolev alternatives.* For $\beta \in \mathbf{N}$ and $1 \leq p < \infty$ with $\beta p > 1$, let

$$\mathcal{F}(\beta, L; p) := \{ l \mid l \text{ is absolutely continuous and } \|l^{(\beta)}\|_p \leq L \},$$

where $\| \cdot \|_p$ denotes the $L_p$-norm. Replacing in the definition of $\rho_n$, $h_n$ and $d_*$ the constant $\beta$ by $\gamma := \beta - 1/p$ and using that $Lh_n^\gamma l(\cdot/h_n) \in \mathcal{F}(\beta, L; p)$ if $l \in \mathcal{F}(\beta, 1; p)$, the results of Theorem 2 extend to Sobolev classes of alternatives as long as the solution of (2) [with a Sobolov ball $\mathcal{F}(\beta, 1; p)$ instead of $\mathcal{H}(\beta, 1)$] has compact support and is of finite total variation. Theorem 3 can be modified in the same way if in addition the corresponding solution of (2) is nonnegative—the final argument in step 3 (proof of Theorem 3) may be replaced with a consideration as in the proof of Theorem 4. The nonnegativity constraint, however, reduces the range of possible Sobolev classes essentially to $\beta = 1$. An explicit solution in case $\beta = 1$ and $p > 2$ has been derived by Sz. Nagy (1941), which satisfies the above requirements in particular.

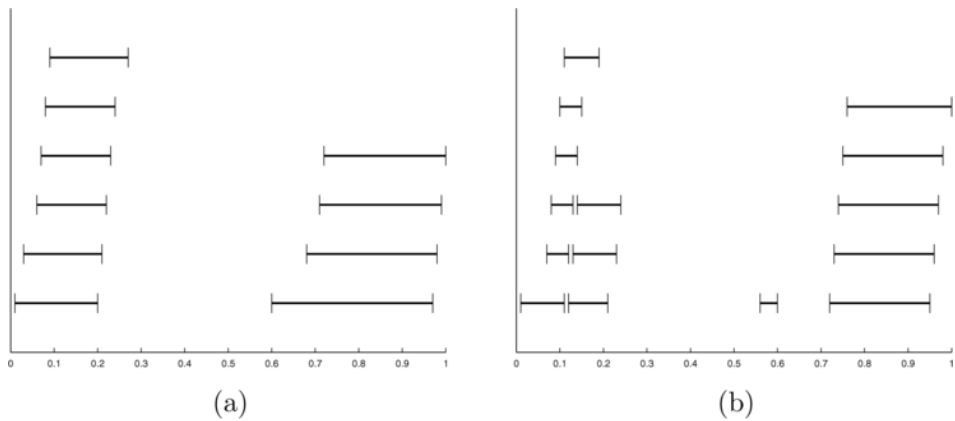

Fig. 2.

5.3. *Random design.* We conjecture that the design assumption (D) can be extended to

(D′) There exists some constant $c > 0$ such that

$$\liminf_{n \to \infty} \frac{H_n(b_n) - H_n(a_n)}{b_n - a_n} \geq c$$

whenever $0 \leq a_n < b_n \leq 1$ and $\liminf_{n \to \infty} \log(b_n - a_n)/\log n > -1$.

Here, $H_n$ denotes the empirical distribution function of the design points. Note that (D) implies (D′). The latter condition is satisfied in particular with probability 1 if $X_1, \ldots, X_n$ are the order statistics of $n$ i.i.d. random variables with a density which is bounded away from zero.

5.4. *Multivariate design.* A further perspective is the extension of the test to two- or even multidimensional design. One application is to detect simultaneously objects on a surface of different shape and size. However, there is no natural class of subsets like intervals one has to look at. Additionally, computational aspects play an increased role: In the univariate case the supremum is taken over $O(n^2)$ single statistics. In two dimensions already, the choice of all rectangles leads to $O(n^4)$.

5.5. *Error laws with point mass and nonsymmetric errors.* If the errors are not restricted to be Lebesgue-continuously distributed, define the local ranks

$$R_{jk}(i) := \sum_{l=j}^{k} \left( 1\{|Y_l| < |Y_i|\} + \frac{1\{|Y_l| = |Y_i|\}}{2} \right) + \frac{1}{2}.$$

The resulting conditional test keeps the significance level.

When the assumption of symmetry is violated, the test is not valid anymore. However, if it seems reasonable in some practical situation that at least $\text{Med}(\varepsilon_i) = 0$, $i = 1, \ldots, n$, one may analyze the data with multiscale sign tests as used in Dümbgen and Johns (2004) for the construction of confidence bands for isotonic median curves. Such a multiscale sign test will be working in a more general setting, but presumably with a considerable loss of efficiency in the Gaussian case.

## 6. Proofs.

PROOF OF THEOREM 1. Let us first introduce some notation. Let $\mathcal{T}_n := \{(j,k) | 1 \leq j < k \leq n\}$ and define the process $X_n$ on $\mathcal{T}_n$ pointwise by

$$X_n(j,k) := \frac{1}{\sqrt{n}} \sum_{i=j}^{k} \psi_{jk}(X_i) \operatorname{sign}(Y_i) \frac{R_{jk}(i)}{k - j + 2}.$$

Since the error distribution is assumed to be symmetric, $\operatorname{sign}(\varepsilon_i)$ is stochastically independent of $|\varepsilon_i|$. Consequently under the null hypothesis, the vector of signs $(\operatorname{sign}(Y_i))_{i=1}^n$ is stochastically independent of the rank vector $R = R_{1n}$. Moreover, $\operatorname{sign}(\varepsilon_i)$ are i.i.d. Rademacher variables. For notational convenience we write $\xi_i$ for $\operatorname{sign}(\varepsilon_i)$.

The proof is partitioned as follows. In step 1, the conditions of Theorem 6.1 in Dümbgen and Spokoiny (2001) are verified for the conditional process $X_n$ given the vector of ranks $R$. Second (step 2), the weak approximation of the conditional process by a Gaussian process in probability is established.

*Step* 1. For any $(j,k) \in \mathcal{T}_n$, let $\sigma_{n,R}^2(j,k)$ denote the conditional variance $\operatorname{Var}(X_n(j,k)|R)$. The sub-Gaussian tails of the conditional process $X_n|R$ are an immediate consequence of Hoeffding's inequality:

$$\mathbb{P}(|X_n(j,k)| > \sigma_{n,R}(j,k)\eta | R)$$
$$= \mathbb{P}\left( \left| \sum_{i=j}^{k} \psi_{jk}(X_i) \xi_i \frac{R_{jk}(i)}{k-j+2} \right| > \left( \sum_{i=j}^{k} \psi_{jk}^2(X_i) \frac{R_{jk}(i)^2}{(k-j+2)^2} \right)^{1/2} \eta \bigg| R \right)$$
$$\leq 2\exp(-\eta^2/2)$$

for any $\eta > 0$, uniformly over $R$ and $1 \leq j < k \leq n$. Let $\rho_n$ be defined by

$$\rho_n((j,k),(j',k'))^2 := |j - j'|/n + |k - k'|/n.$$

In order to show the sub-Gaussian increments of $X_n|R$ with respect to $\rho_n$, it turns out to be sufficient to consider pairs with $j = j' = 1$ and $k < k' = n$, by the same arguments as used in Dümbgen (2002). For any $\eta > 0$, an application of Hoeffding's inequality yields

$$\mathbb{P}\left( \frac{1}{\sqrt{n}} \left| \sum_{i=1}^{n} \psi_{1n}(X_i) \frac{R_{1n}(i)}{n+1} \xi_i - \sum_{i=1}^{k} \psi_{1k}(X_i) \frac{R_{1k}(i)}{k+1} \xi_i \right| > \sqrt{1 - k/n} \eta \bigg| R \right)$$
$$\leq 2\exp(-(1 - k/n)\eta^2/(2B))$$

with

$$B = \operatorname{Var}\left( \frac{1}{\sqrt{n}} \sum_{i=1}^{n} \psi_{1n}(X_i) \frac{R_{1n}(i)}{n+1} \xi_i - \frac{1}{\sqrt{n}} \sum_{i=1}^{k} \psi_{1k}(X_i) \frac{R_{1k}(i)}{k+1} \xi_i \bigg| R \right).$$

First note that $B \leq 2B_1 + 2B_2$, where

$$(3) \quad B_1 = \operatorname{Var}\left(\frac{1}{\sqrt{n}} \sum_{i=1}^{n} \psi_{1n}(X_i) \frac{R_{1n}(i)}{n+1} \xi_i - \frac{1}{\sqrt{n}} \sum_{i=1}^{k} \psi_{1k}(X_i) \frac{R_{1n}(i)}{n+1} \xi_i \bigg| R\right)$$

and

$$(4) \quad B_2 = \operatorname{Var}\left(\frac{1}{\sqrt{n}} \sum_{i=1}^{k} \psi_{1k}(X_i) \frac{R_{1n}(i)}{n+1} \xi_i - \frac{1}{\sqrt{n}} \sum_{i=1}^{k} \psi_{1k}(X_i) \frac{R_{1k}(i)}{k+1} \xi_i \bigg| R\right).$$

Hence it is sufficient to show that $B_i \leq K(1 - k/n)$ for $i = 1, 2$ with some constant $K > 0$ independent of $R$, $k$ and $n$. Throughout this proof, $K$ denotes a generic positive constant depending only on $\psi$ and the design density $h$. Its value may be different in different expressions. Now

$$(5) \quad \begin{aligned} B_1 &= \frac{1}{n} \sum_{i=1}^{k} (\psi_{1n}(X_i) - \psi_{1k}(X_i))^2 \frac{R_{1n}(i)^2}{(n+1)^2} + \frac{1}{n} \sum_{i=k+1}^{n} \psi_{1n}(X_i)^2 \frac{R_{1n}(i)^2}{(n+1)^2} \\ &\leq \frac{1}{n} \sum_{i=1}^{k} (\psi_{1n}(X_i) - \psi_{1k}(X_i))^2 + K(1 - k/n). \end{aligned}$$

For notational convenience, we denote the scale $(X_k - X_1)$ by $t_{1k}$. The finite total variation of $\psi$ implies that $\psi(x) = \int_{[0,x]} g(u)\,dP(u)$ for all but at most countably many numbers $x \in [0,1]$, where $P$ is some probability measure on $[0,1]$ and $g$ is some measurable function with $|g| \leq TV(\psi)$. For $0 \leq z_1 \leq z_2 \leq 1$ let $\mu$ be defined by $\mu([z_1, z_2]) := \int_{z_1}^{z_2} |g(x)| P(dx)$. Note that $|\psi(z_1) - \psi(z_2)| \leq \mu([z_1, z_2])$. Let $H_n$ denote the empirical distribution function of the design points and define

$$A_x^{(kn)} := \left[\frac{x - X_1}{t_{1n}}, \frac{x - X_1}{t_{1k}}\right].$$

The sum in (5) is then bounded by

$$(6) \quad \frac{1}{n} \sum_{i=1}^{k} (\psi_{1n}(X_i) - \psi_{1k}(X_i))^2$$

$$= \frac{1}{n} \sum_{i=1}^{k} \{\psi((X_i - X_1)/t_{1n}) - \psi((X_i - X_1)/t_{1k})\}^2$$

$$\leq \int_{X_1}^{X_k} \mu(A_x^{(kn)})^2 H_n(dx)$$

$$(7) \quad = \int \left\{\int_{X_1}^{X_k} I\{y \in A_x^{(kn)}, z \in A_x^{(kn)}\} H_n(dx)\right\} \mu(dy)\mu(dz)$$

$$\leq K \sup_{y\in[0,1]} \int_{X_1}^{X_k} I\{y \in A_x^{(kn)}\} H_n(dx)$$

(8) $$\leq K \sup_{y\in[0,1]} (H_n(yt_{1n} + X_1) - H_n(yt_{1k} + X_1)),$$

where equality (7) follows by an application of Fubini's theorem. But the design assumption (D) implies that $H - 1/n \leq H_n \leq H$ pointwise. Therefore, the latter supremum in (8) is bounded by

$$\sup_{y\in[0,1]} (H(yt_{1n} + X_1) - H(yt_{1k} + X_1)) + 1/n \leq K \int_{X_k}^{X_n} h(x)\lambda(dx) + 1/n,$$

which is bounded from above by $K(1 - k/n)$ for some constant $K$ independent of $n$ and $k$. In order to bound $B_2$ in (4), define $\tilde{R}_{1k}(i) := \sum_{l=k+1}^{n} I\{|Y_l| \leq |Y_i|\}$; thus $R_{1n}(i)$ equals $R_{1k}(i) + \tilde{R}_{1k}(i)$ a.s. Then

$$B_2 \leq \frac{2}{n} \sum_{i=1}^{k} \psi_{1k}(X_i)^2 \left(\frac{k+1}{n+1} - 1\right)^2 \frac{R_{1k}(i)^2}{(k+1)^2} + \frac{2}{n} \sum_{i=1}^{k} \psi_{1k}(X_i)^2 \frac{\tilde{R}_{1k}(i)^2}{(n+1)^2}$$

$$\leq K(1 - k/n)^2 + K\frac{2}{n} \sum_{i=k+1}^{n} \frac{i^2}{(n+1)^2}$$

$$\leq K(1 - k/n).$$

Consequently, $X_n|R$ has sub-Gaussian increments with respect to $\rho_n$.

For some totally bounded pseudometric space $(\mathcal{T}, \rho)$, $\mathcal{T}' \subset \mathcal{T}$ and any $\varepsilon > 0$, the covering number $N(\varepsilon, \mathcal{T}', \rho)$ is defined as the infimum of $\sharp \mathcal{T}_0$ over all $\mathcal{T}_0 \subset \mathcal{T}'$ such that $\inf_{t_0 \in \mathcal{T}_0} \rho(t_0, t) \leq \varepsilon \ \forall t \in \mathcal{T}'$. To finish step 1, we need to establish the bound for the covering numbers,

$$N((\delta u)^{1/2}, \{(j,k) \in \mathcal{T}_n : \sigma(j,k)_{n,R}^2 \leq \delta\}, \rho_n) \leq A u^{-2} \delta^{-1}$$

with a constant $A > 0$, independent of $R$ and $n$. Since $\psi$ is continuous with $\int_0^1 \psi(x)\,dx > 0$, there exists some nondegenerate interval $[a,b] \subset [0,1]$ with $\psi(x)^2 \geq \tau$ for some strictly positive constant $\tau$ and any $x \in [a,b]$. Let $B_{jk} := \{i : (X_i - X_j)/t_{jk} \in [a,b]\}$. By assumption (D),

$$\frac{\sharp B_{jk}}{n} = \int_{t_{jk}a + X_j}^{t_{jk}b + X_j} dH_n(x)$$

$$\geq H(t_{jk}b + X_j) - H(t_{jk}a + X_j) - \frac{1}{n} \geq K\frac{k - j - 1/K}{n}.$$

This entails the lower bound

$$\sigma_{n,R}(j,k)^2 \geq \frac{1}{n} \sum_{i \in B_{jk}} \tau \frac{R_{jk}(i)^2}{(k-j+2)^2}$$

$$\geq \frac{\tau}{n} \sum_{i=1}^{\sharp B_{jk}} \frac{i^2}{(k-j+2)^2}$$

$$= \frac{1}{n} \tau \frac{(\sharp B_{jk})(\sharp B_{jk}+1)(2\sharp B_{jk}+1)}{6(k-j+2)^2} \geq K \frac{k-j-1/K}{n},$$

with some constant $K > 0$, independent of $R, k, j$ and $n$. Therefore,

$$N((\delta u)^{1/2}, \{(j,k) \in \mathcal{T}_n : \sigma_{n,R}(j,k)^2 \leq \delta\}, \rho_n)$$
$$\leq N((\delta u)^{1/2}, \{(j,k) \in \mathcal{T}_n : (k-j)/n \leq (\delta+1/n)/K\}, \rho_n).$$

If $\delta \geq 1/n$, then $\delta + 1/n \leq 2\delta$, and via the embedding $k \mapsto k/n$ of $\mathcal{T}_n$ into $[0,1]$, the covering number can be bounded by $Au^{-2}\delta^{-1}$ for some constant $A > 0$ with the same argument as given in Dümbgen and Spokoiny (2001). Note that the desired bound is necessarily satisfied for $\delta \leq 1/n$: Then $\sharp\{(j,k) \in \mathcal{T}_n : (k-j)/n \leq (\delta+1/n)/K\} \leq \sharp\{(j,k) \in \mathcal{T}_n : (k-j) \leq 2/K\} \leq 2K^{-1}n \leq 2K^{-1}\delta^{-1}$.

*Step* 2. Let $\mathcal{S}_n := \{(X_i, X_j) | 0 \leq j < k \leq n\}$, where $X_0 := 0$. Redefine the process $X_n$ on $\mathcal{S}_n$ via

$$X_n(s,t) := \frac{1}{\sqrt{n}} \sum_{i \in I_{st}} \psi_{st}(X_i) \xi_i \frac{R_{st}(i)}{\sharp I_{st}+1}, \qquad (s,t) \in \mathcal{S}_n,$$

where $I_{st} := \{i | X_i \in [s,t]\}$ and $R_{st}$ denotes the rank of $|Y_i|$ among the $\sharp I_{st}$ numbers $|Y_k| : X_k \in [s,t]$. Furthermore, let the process $Z$ on $\mathcal{S} := \{(s,t) | 0 \leq s < t \leq 1\}$ pointwise be defined by

$$Z(s,t) := \frac{1}{\sqrt{3}} \int_s^t \psi_{st}(x) \sqrt{h(x)} \, dW(x), \qquad (s,t) \in \mathcal{S},$$

with $W$ some Brownian motion on the unit interval. In the sequel we prove the weak convergence in probability of the conditional process under the null hypothesis, that is,

$$d_w(\mathcal{L}(X_n|R), \mathcal{L}(Z(s,t))_{(s,t) \in \mathcal{S}_n}) \longrightarrow_p 0,$$

where $d_w$ denotes some metric generating the topology of weak convergence. It follows by a standard chaining argument and the above established results that uniformly over $R$ and $n$, $X_n | R$ is stochastically equicontinuous with respect to $\rho$, pointwise defined by

$$\rho((s,t),(s',t'))^2 := |H(s) - H(s')| + |H(t) - H(t')|.$$

To prove the weak convergence in probability, it is therefore sufficient to show the convergence of the finite-dimensional distributions of $X_n | R$. Let

$$\phi_{i,n}(s,t) := \frac{1}{\sqrt{n}} I_{[s,t]}(X_i) \psi_{st}(X_i) \xi_i \frac{R_{st}(i)}{\sharp I_{st}+1}, \qquad (s,t) \in \mathcal{S}_n.$$

Then $X_n(s,t) = \sum_{i=1}^n \phi_{i,n}(s,t)$, and the $\phi_{i,n}$ are independent conditioned on $R$. One verifies that

$$\mathbb{E}\left(\sum_{i=1}^n \|\phi_{i,n}\|_{\mathcal{S}_n}^2 \Big| R\right) \leq \|\psi\|_{\sup}^2$$

and for arbitrary $u > 0$,

$$\mathbb{E}\left(\sum_{i=1}^n I\{\|\phi_{i,n}\|_{\mathcal{S}_n}^2 > u\}\|\phi_{i,n}\|_{\mathcal{S}_n}^2 \Big| R\right) = o(1).$$

For any natural number $k$, let now $\{(s_1,t_1),\ldots,(s_k,t_k)|0 \leq s_i < t_i \leq 1, i = 1,\ldots,k\}$ and $\mathcal{S}_n^k = \{(s_{1n},t_{1n}),\ldots,(s_{kn},t_{kn})\} \subset \mathcal{S}_n$ such that $(s_{ni},t_{ni}) \to (s_i,t_i)$ for $i=1,\ldots,k$. For a given vector $R$ of ranks, let us introduce the process $Z_{nR}$ on $\mathcal{S}_n$ which is, conditioned on $R$, a centered Gaussian process with conditional covariance structure as $X_n|R$, that is,

(9)
$$\operatorname{cov}(Z_{nR}(s,t), Z_{nR}(s',t')|R)$$
$$= \frac{1}{n}\sum_{i \in I_{st} \cap I_{s't'}} \psi_{st}(X_i)\psi_{s't'}(X_i) \frac{R_{st}(i)}{\sharp I_{st}+1} \frac{R_{s't'}(i)}{\sharp I_{s't'}+1}.$$

Since the conditional covariance function of $X_n|R$ is uniformly bounded by $\pm\|\psi\|_{\sup}^2$, respectively, Lindeberg's central limit theorem entails that

$$d_w(\mathcal{L}(X_{n|\mathcal{S}_n^k}|R), \mathcal{L}(Z_{nR|\mathcal{S}_n^k}|R)) \to 0,$$

due to the compactness of $[-\|\psi\|_{\sup}^2, \|\psi\|_{\sup}^2]$. It remains to be shown that

(10) $$d_w(\mathcal{L}(Z_{nR|\mathcal{S}_n^k}|R), \mathcal{L}(Z_{n|\mathcal{S}_n^k})) \longrightarrow_p 0.$$

Let $(s_n,t_n) \in \mathcal{S}_n$ with $\liminf_n |s_n - t_n| > 0$. Then

$$\left|\frac{R_{s_n t_n}(i)}{\sharp I_{s_n t_n}+1} - (F(|Y_i|) - F(-|Y_i|))\right|$$
$$\leq \sup_z \left|\frac{1}{\sharp I_{s_n t_n}+1}\sum_{j \in I_{s_n t_n}} I\{|Y_j| \leq |z|\} - (F(|z|) - F(-|z|))\right|$$

and the latter quantity is $o_p(1)$ by the Glivenko–Cantelli theorem. This shows that for $(s_n,t_n),(s'_n,t'_n) \in \mathcal{S}_n$ with $(s_n,t_n) \to (s,t) \in \mathcal{S}$ and $(s'_n,t'_n) \to (s',t') \in \mathcal{S}$, (9) is equal to

$$\operatorname{cov}(X_n(s_n,t_n), X_n(s'_n,t'_n)|R)$$
$$= n^{-1}\sum_{i \in I_{s_n t_n} \cap I_{s'_n t'_n}} \psi_{s_n t_n}(X_i)\psi_{s'_n t'_n}(X_i)(F(|Y_i|) - F(-|Y_i|))^2 + o_p(1).$$

The random variables $\operatorname{sign}(Y_i)\{F(|Y_i|) - F(-|Y_i|)\} = 2F(Y_i) - 1$, $i = 1, \ldots, n$, are independent and uniformly distributed on $[-1, 1]$. Consequently, assumption (D) and an application of Chebyshev's inequality finally yields

$$\operatorname{cov}(X_n(s_n, t_n), X_n(s'_n, t'_n)|R) \longrightarrow_p \tfrac{1}{3} \int \psi_{st}(x) \psi_{s't'}(x) h(x)\, dx$$

which implies (10).

From steps 1 and 2 the asserted stochastically weak convergence of our test statistic can be deduced with the same argument as given in Dümbgen (2002), page 528. □

PROOF OF THEOREM 2. For a fixed smoothness class $\mathcal{H}(\beta, L)$, let $\gamma = \gamma_\beta$ be the solution of the optimization problem (2). As pointed out in Section 2, $\gamma$ is an even function with compact support, say $[-C, C]$. Now define the following set of testing functions: For a given bandwidth $h_n > 0$ and any integer $j$ let

$$\gamma_{j,n}(\cdot) := \gamma\left(\frac{\cdot - (2j-1)Ch_n}{h_n}\right) \quad \text{and define} \quad g_{j,n}(\cdot) := \frac{1}{\sqrt{h(\cdot)}} L h_n^\beta \gamma_{j,n}.$$

[Note that $h(\cdot)$ denotes the design density whereas $h_n$ denotes the $n$-dependent scale parameter.] Let $[a, a+b] \subset J$ for some $b > 0$ and define

$$\mathcal{J}_n := \{j \in \mathbb{N} : (2j-1)Ch_n \in [a + Ch_n, a + b - Ch_n]\}.$$

Let $\mathcal{G}_n := \{g_{j,n} : j \in \mathcal{J}_n\}$. Note that $g \in \mathcal{H}_h(\beta, L)$ for every $g \in \mathcal{G}_n$. Following the arguments in Dümbgen and Spokoiny (2001), proof of Theorem 3.1a, one shows that for any test $\phi : \mathbb{R}^n \to [0, 1]$ with significance level $\leq \alpha$,

$$\inf_{g \in \mathcal{G}_n} \mathbb{E}_g \phi(X, Y) - \alpha \leq \mathbb{E}_0 \left| \frac{1}{\sharp \mathcal{G}_n} \sum_{g \in \mathcal{G}_n} \frac{d\mathbb{P}_g}{d\mathbb{P}_0}(X, Y) - 1 \right|.$$

The aim is to determine $h_n$ such that the right-hand side tends to zero as $n$ goes to infinity. Define the index set $I_g := \{i | g(X_i) > 0\}$. By construction, $I_g \cap I_{g'} = \varnothing$ for $g \neq g'$ and $g, g' \in \mathcal{G}_n$. Then for any $g \in \mathcal{G}_n$, the likelihood ratio equals

$$\frac{d\mathbb{P}_g}{d\mathbb{P}_0}(X, Y) = \prod_{i \in I_g} \frac{f(Y_i - g(X_i))}{f(Y_i)},$$

which shows that $d\mathbb{P}_g/d\mathbb{P}_0(X, Y)$, $g \in \mathcal{G}_n$, are independent. Note that their expectation is not the same for every $g$. Using a standard truncation argument as in Dümbgen and Walther (2008), proof of Lemma 10, it turns out to be sufficient to find $h_n$ such that

$$
\begin{aligned}
(11) \quad &\inf_{\delta \in (0,\delta_0]} \max_{g \in \mathcal{G}_n} \frac{1}{(\sharp \mathcal{G}_n)^\delta} \mathbb{E}_0\left(\left(\frac{d\mathbb{P}_g}{d\mathbb{P}_0}(X,Y)\right)^{1+\delta}\right) \\
&= \inf_{\delta \in (0,\delta_0]} \max_{g \in \mathcal{G}_n} \frac{1}{(\sharp \mathcal{G}_n)^\delta} \prod_{i=1}^n \left\{\int f(y)\left(\frac{f(y-g(X_i))}{f(y)}\right)^{1+\delta} dy\right\} \to 0
\end{aligned}
$$

as $n \to \infty$. Using the expansion of assumption (E2), (11) is equal to

$$
\inf_{\delta \in (0,\delta_0]} \max_{g \in \mathcal{G}_n} \frac{1}{(\sharp \mathcal{G}_n)^\delta} \prod_{i=1}^n \left\{1 + \frac{1}{2}\delta(1+\delta)I(f)g(X_i)^2(1+r(g(X_i),\delta))\right\}.
$$

But for $h_n$ sufficiently small, the latter expression is bounded by

$$
(12) \quad \inf_{\delta \in (0,\delta_0]} \max_{g \in \mathcal{G}_n} \exp(n\tfrac{1}{2}\delta(1+\delta)I(f)\|g\|_{n,2}^2(1+\bar{r}(g)) - \delta\log(\sharp \mathcal{G}_n)),
$$

using the series representation of the logarithm, where

$$
\|g\|_{n,2} := \frac{1}{n}\sum_{i=1}^n g(X_i)^2
$$

and $\bar{r}(g) := \sup_{\delta \in (0,\delta_0]} \sup_{x \in [0,1]} |r(g(x),\delta)|$. Furthermore,

$$
\begin{aligned}
&\frac{1}{n}\sum_{i=1}^n g_{j,n}(X_i)^2 - \int g_{j,n}(x)^2 h(x)\, dx \\
&= L^2 h_n^{2\beta} \sum_{i \in I_{g_{j,n}}} \int_{X_{i-1}}^{X_i} \left(\frac{\gamma_{j,n}(X_i)^2}{h(X_i)} - \frac{\gamma_{j,n}(x)^2}{h(x)}\right)h(x)\, dx \\
&\leq L^2 h_n^{2\beta} \sum_{i \in I_{g_{j,n}}} \sup_{x \in [X_{i-1},X_i]} \left|\frac{\gamma_{j,n}(X_i)^2}{h(X_i)} - \frac{\gamma_{j,n}(x)^2}{h(x)}\right| \frac{1}{n}.
\end{aligned}
$$

The last expression is of order $O(h_n^{2\beta}n^{-1})$: Since the design density $h$ is of bounded total variation as well as uniformly bounded away from zero, also $1/h$ is of bounded total variation. In addition, $\gamma$ is bounded and of bounded total variation (for $\beta \leq 1$, $\gamma$ is explicitly known and unimodal, while its first derivative is Hölder-continuous in case $\beta > 1$). Consequently, $TV(\gamma_{j,n}^2/h) \leq K(TV(\gamma^2) + TV(h)) < \infty$ with some constant $K$ independent of $j$ and $n$, which shows that $\|g_{j,n}\|_{n,2}^2 = h_n^{2\beta+1}\|\gamma\|_2^2(1+O((h_n n)^{-1}))$. Thus (12) is bounded by

$$
(13) \quad \inf_{\delta \in (0,\delta_0]} \max_{g \in \mathcal{G}_n} \exp(n\tfrac{1}{2}\delta(1+\delta)I(f)L^2 h_n^{2\beta+1}\|\gamma\|_2^2(1+R(n,g))
$$
$$
- \delta\log(\sharp \mathcal{G}_n)),
$$

with a sequence $R(n,g)$ of order $O(\max\{(h_n n)^{-1}, \bar{r}(g)\})$.

Let $\varepsilon_n > 0$ be arbitrary numbers with $\varepsilon_n \to 0$ and $\varepsilon_n \sqrt{\log n} \to \infty$. Define the bandwidth

$$h_n := \left(\frac{d_* \rho_n}{L}\right)^{1/\beta} (1 - \varepsilon_n)^{1/\beta},$$

which implies that $\sup_{g \in \mathcal{G}_n} R(n,g)$ in (13) is of order $(\log n)^{-1}$. By the choice of $\mathcal{G}_n$, $\sharp\mathcal{G}_n \geq b/(2Ch_n) - 1$. Let $\delta = \delta_n := \varepsilon_n$. Then (13) is bounded by

$$\exp(\varepsilon_n(1+\varepsilon_n)(2\beta+1)^{-1}\log n(1-\varepsilon_n)^{(2\beta+1)/\beta}$$
$$\qquad - \varepsilon_n(2\beta+1)^{-1}(\log n - \log\log n) + o(1))$$
$$= \exp\left(-\frac{1+\beta}{\beta}\varepsilon_n^2(1+O(\varepsilon_n))\log n + \varepsilon_n(2\beta+1)^{-1}\log\log n + o(1)\right),$$

which tends to zero as $n$ goes to infinity. □

PROOF OF THEOREM 3. By virtue of the proof of Theorem 1, the conditional process $X_n|R$ satisfies the conditions of Theorem 6.1 of Dümbgen and Spokoiny (2001) uniformly in $R$ and $n$. This entails that there exists some constant $C > 0$ independent of $n$ with $\kappa_\alpha^n(R) \leq C$, where $\kappa_\alpha^n(R)$ denotes the $(1-\alpha)$-quantile of $\mathcal{L}(T_n|R)$ under the null hypothesis. Consequently,

$$\mathbb{P}_l(\phi_n^* = 1) = \int \mathbb{P}_l(T_n > \kappa_\alpha^n(R)|R)\, d\mathbb{P}_l(R)$$
$$\geq \int \mathbb{P}_l(T_n > C|R)\, d\mathbb{P}_l(R) = \mathbb{P}_l(T_n > C).$$

Furthermore, $\mathbb{P}_l(T_n > C) \geq \mathbb{P}_l(|T_{jk}| > C + \sqrt{2\log(n/(k-j))})$ for any $1 \leq j < k \leq n$. It is therefore sufficient to show that for any sequence $l_n \in \mathcal{H}_h(\beta, L)$ with maximal absolute value $\|l_n\sqrt{h}\|_{\sup} \geq d^*\rho_n(1+\varepsilon_n)$, there exists a sequence of pairs $(j_n, k_n)$ with $1 \leq j_n < k_n \leq n$ such that

$$\liminf_{n\to\infty} \mathbb{P}_{l_n}(|T_{j_n k_n}| > C + \sqrt{2\log(n/(k_n - j_n))}) = 1.$$

The proof is organized as follows: At first (step 1), the $L_2$-approximation of the numerator of $T_{j_n k_n}$ by a sum of independent random variables is established. Second (step 2), Taylor-type expansions of its expectation and variance are provided, and the asymptotic power of our test is determined along sequences of alternatives converging to zero at the fastest possible rate. Finally (step 3), we treat alternatives converging to zero at a slow rate or staying uniformly bounded away from zero.

*Step* 1. Let $I_n := \{j_n, \ldots, k_n\}$ be an interval of indices with $1 \leq j_n < k_n \leq n$ and $\sharp I_n = k_n - j_n + 1 \to \infty$. For notational convenience, denote $\psi_n := \psi_{j_n k_n}$ and $R_n(i) := R_{j_n k_n}(i), i \in I_n$. Let $S_n$ be the (normalized) numerator of the single local test statistic $T_{j_n k_n}$, that is,

$$S_n := \frac{1}{\sqrt{\sharp I_n}} \sum_{i \in I_n} \psi_n(X_i) \operatorname{sign}(Y_i) \frac{R_n(i)}{\sharp I_n + 1}$$

(14)

$$\stackrel{a.s.}{=} \sum_{i \in I_n} \sum_{j \in I_n} \frac{1}{\sharp I_n + 1} \operatorname{sign}(Y_i) \frac{\psi_n(X_i)}{\sqrt{\sharp I_n}} I\{|Y_j| \leq |Y_i|\}.$$

In the sequel, we establish the approximation of $S_n$ by a sum of independent random variables which is up to $O_p(1/\sharp I_n)$ its Hájek projection [see, e.g., van der Vaart (1998)]. For that purpose the Hoeffding decomposition is applied. With $c_i = c_{n,i} := (\sharp I_n)^{-1/2} \psi_n(X_i)$, let $A_{ij} := \operatorname{sign}(Y_i) c_i I\{|Y_j| \leq |Y_i|\}$ and define $H_{ij} := A_{ij} + A_{ji}$. Then

$$S_n \stackrel{a.s.}{=} \sum_{i \in I_n} \sum_{\substack{j \in I_n: \\ j < i}} \frac{1}{\sharp I_n + 1} H_{ij} + \sum_{i \in I_n} \frac{1}{\sharp I_n + 1} A_{ii}.$$

With the definition

$$\tilde{H}_{ij} := \mathbb{E}(S_n | Y_i, Y_j) - \mathbb{E}(S_n | Y_i) - \mathbb{E}(S_n | Y_j) + \mathbb{E}(S_n)$$
$$= H_{ij} - \mathbb{E}(H_{ij} | Y_i) - \mathbb{E}(H_{ij} | Y_j) + \mathbb{E} H_{ij}$$

for $i \neq j$, we obtain the decomposition

$$S_n \stackrel{a.s.}{=} \sum_{i \in I_n} \sum_{\substack{j \in I_n: \\ j < i}} \frac{1}{\sharp I_n + 1} \tilde{H}_{ij}$$
$$+ \sum_{i \in I_n} \left( \frac{H_{ii}/2}{\sharp I_n + 1} + \sum_{\substack{j \in I_n: \\ j < i}} \frac{1}{\sharp I_n + 1} (\mathbb{E}(H_{ij}|Y_i) + \mathbb{E}(H_{ij}|Y_j) - \mathbb{E} H_{ij}) \right)$$
$$=: S_n^{(0)} + S_n^{(1)},$$

where $S_n^{(0)}$ and $S_n^{(1)}$ are uncorrelated. Note that in particular $\mathbb{E} \tilde{H}_{ij} = 0$ and $\operatorname{cov}(\tilde{H}_{ij}, \tilde{H}_{kl}) = 0$ for $(i,j) \neq (k,l)$. Consequently

$$\operatorname{Var}(S_n - S_n^{(1)}) = \frac{1}{(\sharp I_n + 1)^2} \sum_{i \in I_n} \sum_{j \in I_n: j < i} \operatorname{Var}(\tilde{H}_{ij})$$
$$\leq \frac{1}{(\sharp I_n + 1)^2} \sum_{i \in I_n} \sum_{j \in I_n: j < i} 4 c_{n,i}^2$$
$$= O(1/\sharp I_n),$$

since by construction, $\mathrm{Var}(\tilde{H}_{ij}) \leq \mathrm{Var}(H_{ij})$. Furthermore, $S_n^{(1)}$ is equal to

$$\sum_{i \in I_n} \frac{c_i}{\sharp I_n + 1} \mathrm{sign}(Y_i)$$

$$+ \sum_{\substack{i,j \in I_n \\ j \neq i}} \frac{c_i}{\sharp I_n + 1} \mathrm{sign}(Y_i)(F_j(|Y_i|) - F_j(-|Y_i|))$$

$$+ \sum_{\substack{i,j \in I_n \\ j \neq i}} \frac{c_i}{\sharp I_n + 1} \left\{ \int_{\mathbb{R} \setminus [-|Y_j|,|Y_j|]} \mathrm{sign}(y) \, dF_i(y) - \mathbb{E}(H_{ij}) \right\},$$

where $F_i$ denotes the distribution function of $Y_i$. For any distribution function $F$, let $G$ be pointwise defined on $\mathbb{R}^+$ by $G(t) := F(t) - F(-t-)$, with $F(y-)$ the limit on the left, that is, $\lim_{x \nearrow y} F(x)$. We denote $\bar{F} := 1/(\sharp I_n) \sum_{i \in I_n} F_i$, $\bar{G}(t) := \bar{F}(t) - \bar{F}(-t-)$ and $\bar{F}^\psi := 1/(\sharp I_n) \sum_{i \in I_n} \psi_n(X_i) F_i$. Then $\mathbb{E}(S_n^{(1)} - \widehat{S}_n)^2 = O(1/\sharp I_n)$, with

$$\widehat{S}_n := \frac{1}{\sqrt{\sharp I_n}} \sum_{i \in I_n} \left\{ \psi_n(X_i) \mathrm{sign}(Y_i) \bar{G}(|Y_i|) \right.$$

(15)
$$+ \int_{\mathbb{R} \setminus [-|Y_i|,|Y_i|]} \mathrm{sign}(y) \, d\bar{F}^\psi(y)$$

$$\left. - \mathbb{E} \int_{\mathbb{R} \setminus [-|Y_i|,|Y_i|]} \mathrm{sign}(y) \, d\bar{F}^\psi(y) \right\}.$$

*Step* 2. For two functions $f$ and $g$ in $L_2[0,1]$, let

$$\langle f, g \rangle_{I_n} := 1/(\sharp I_n) \sum_{i \in I_n} f(X_i) g(X_i)$$

and let $\|f\|_{I_n,2} := \langle f, f \rangle_{I_n}^{1/2}$ denote the corresponding norm. Let $(l_n)$ be a sequence of alternatives. If $M(l_n)$ denotes the maximal point of $|l_n|$, let $(X_{j_n}, X_{k_n})$ be the design points which are closest to $M(l_n) - h_n$ and $M(l_n) + h_n$, respectively, where $h_n := (\delta_n/L)^{1/\beta}$ with $\delta_n := d^* \rho_n (1 + \varepsilon_n)$. Symmetry considerations show that we may assume without loss of generality that $l_n$ is positive at $M(l_n)$. Besides the restriction $\|l_n \sqrt{h}\|_{\sup} \geq d^* \rho_n (1 + \varepsilon_n)$, it is assumed in this paragraph that

(16) $$\|l_n\|_{\sup}/\rho_n = O(1),$$

which is equivalent to $\|l_n \sqrt{h}\|_{\sup}/\rho_n = O(1)$. Note that (16) implies $\sqrt{\sharp I_n} \times \|l_n\|_{I_n,2}^2 = o(1)$.

Our first goal is to show that

$$\text{(17)} \quad \frac{\mathbb{E}_{l_n}\widehat{S}_n}{\sqrt{\text{Var}_{l_n}\widehat{S}_n}} = \sqrt{12}\sqrt{\sharp I_n}\frac{\langle \psi_n, l_n\rangle_{I_n}}{\|\psi_n\|_{I_n,2}}\int f(y)^2\,dy + o(1)$$

for any sequence $(l_n)$ satisfying (16). The symmetry of the error distribution around zero and the boundedness of the first derivative $f'$ provide the expansion

$$\text{sign}(Y_i)\bar{G}(|Y_i|)$$

$$= \text{sign}(Y_i)\bigg\{(F(|Y_i|) - F(-|Y_i|))$$

$$- (f(|Y_i|) - f(-|Y_i|))\bigg(\frac{1}{\sharp I_n}\sum_{j\in I_n}l_n(X_j)\bigg) + O_{\text{unif}}(\|l_n\|_{I_n,2}^2)\bigg\}$$

$$= (2F(Y_i) - 1) + O_{\text{unif}}(\|l_n\|_{I_n,2}^2).$$

Here and in what follows, a sequence of random variables $(Z_n)$ is $O_{\text{unif}}(c_n)$ with a sequence of positive numbers $(c_n)$, if $\limsup_n |Z_n/c_n| \leq c < \infty$ with some nonrandom nonnegative constant $c$. In order to treat the expectation

$$\mathbb{E}_{l_n}\widehat{S}_n = \frac{1}{\sqrt{\sharp I_n}}\sum_{i\in I_n}\psi_n(X_i)\bigg\{\int (2F(y) - 1)\,dF_i(y) + O(\|l_n\|_{I_n,2}^2)\bigg\},$$

first observe that for any $\theta \in \mathbb{R}$, $\int_{\mathbb{R}}(2F(y) - 1)f(y + \theta)\,dy = \int_{\mathbb{R}}f'(t) \times \int_{t-\theta}^{t}(2F(y) - 1)\,dy\,dt$, using Fubini's theorem and the symmetry of the error density $f$. Taylor expansion of the inner integral entails that

$$\mathbb{E}_{l_n}\widehat{S}_n = \sqrt{\sharp I_n}\langle \psi_n, l_n\rangle_{I_n}\bigg\{-\int (2F(y) - 1)f'(y)\,dy\bigg\} + \sqrt{\sharp I_n}O(\|l_n\|_{I_n,2}^2)$$

$$\text{(18)}$$

$$= 2\sqrt{\sharp I_n}\langle \psi_n, l_n\rangle_{I_n}\bigg\{\int f(y)^2\,dy\bigg\} + \sqrt{\sharp I_n}O(\|l_n\|_{I_n,2}^2),$$

where the last equality is obtained via partial integration. Furthermore,

$$\text{Var}_{l_n}\bigg(\frac{1}{\sqrt{\sharp I_n}}\sum_{i\in I_n}\psi_n(X_i)\text{sign}(Y_i)\bar{G}(|Y_i|)\bigg)$$

$$\text{(19)}$$

$$= \frac{1}{\sharp I_n}\sum_{i\in I_n}\psi_n(X_i)^2\mathbb{E}_{l_n}(2F(Y_i) - 1)^2 + O(\|l_n\|_{I_n,2}^2).$$

In order to bound the variance of the second part in the approximation (15), namely

$$\text{(20)} \quad \text{Var}_{l_n}\bigg(\frac{1}{\sqrt{\sharp I_n}}\sum_{i\in I_n}\int_{\mathbb{R}\setminus[-|Y_i|,|Y_i|]}\text{sign}(y)\,d\bar{F}^\psi(y)\bigg)$$

$$\text{(21)} \qquad \leq \frac{1}{\sharp I_n} \sum_{i \in I_n} \mathbb{E}_{l_n} \left( \int_{\mathbb{R} \setminus [-|Y_i|, |Y_i|]} \text{sign}(y) \, d\bar{F}^\psi(y) \right)^2,$$

note that by the symmetry of $\text{sign}(\cdot)$ and Fubini's theorem,

$$\left| \int_{[-z,z]^c} \text{sign}(y) \, d\bar{F}^\psi(y) \right|$$
$$= \left| \frac{1}{\sharp I_n} \sum_{i \in I_n} \psi_n(X_i) \int_{\mathbb{R}} f'(t) \int_{[-z,z]^c} -\text{sign}(y) I\{y \in [t, t + l_n(X_i)]\} \, dy \, dt \right|$$
$$\leq \langle \psi_n, |l_n| \rangle_{I_n} \int_{\mathbb{R}} |f'(t)| \, dt.$$

This shows that (20) is $O(\|l_n\|_{I_n,2}^2)$ by Cauchy–Schwarz. Furthermore,

$$\int_{\mathbb{R}} (2F(y) - 1)^2 \, d(F_i(y) - F(y))$$
$$= \int_{\mathbb{R}} (2F(y) - 1)^2 \int_{y - l_n(X_i)}^{y} -f'(t) \, dt \, dy$$
$$= \int_{\mathbb{R}} f'(t) \int_{t}^{t + l_n(X_i)} -(2F(y) - 1)^2 \, dy \, dt$$
$$= l_n(X_i) \int_{\mathbb{R}} 4f(t)^2 (2F(t) - 1) \, dt + O(l_n(X_i)^2),$$

where the latter integral is equal to zero by the symmetry of the error distribution. This finally gives together with (19) and the bound of (20)

$$\text{(22)} \qquad \text{Var}_{l_n} \widehat{S}_n = \tfrac{4}{12} \|\psi_n\|_{I_n,2}^2 + O(\|l_n\|_{I_n,2}^2).$$

Note at this point that $\text{Var}_{l_n} \widehat{S}_n$ is uniformly bounded from above and from below. Thus the combination of (18) and (22) entails (17) for any sequence $(l_n)$ satisfying (16).

In the next step, it will be shown that the denominator of $T_{j_n k_n}$ is a sufficiently good approximation for the standard deviation of $\widehat{S}_n$ under the sequence of alternatives $l_n$. Remember that it is the conditional standard deviation given the vector of ranks of the numerator under the null hypothesis. Using the representation $R_n(i) = \sum_{k \in I_n} I\{|Y_k| \leq |Y_i|\}$ a.s., one verifies that

$$\mathbb{E}_{l_n} \left( \frac{1}{\sharp I_n} \sum_{i \in I_n} \psi_n(X_i)^2 \frac{R_n(i)^2}{(\sharp I_n + 1)^2} \right) = \frac{4}{12} \|\psi_n\|_{I_n,2}^2 + O(\|l_n\|_{I_n,2}^2),$$

and analogously for $i, j \in I_n$ with $i \neq j$

$$\mathbb{E}_{l_n} \left( \frac{R_n(i)^2}{(\sharp I_n + 1)^2} \frac{R_n(j)^2}{(\sharp I_n + 1)^2} \bigg| Y_i, Y_j \right) = \bar{G}(|Y_i|)^2 \bar{G}(|Y_j|)^2 + O_{\text{unif}}(1/\sharp I_n)$$

and

$$\mathrm{Var}_{l_n}\left(\frac{1}{\sharp I_n}\sum_{i\in I_n}\psi_n(X_i)^2\frac{R_n(i)^2}{(\sharp I_n+1)^2}\right)=O(1/\sharp I_n),$$

which by Chebyshev's inequality shows in particular that under condition (16)

(23)
$$\sqrt{\sharp I_n}\frac{\langle\psi_n,l_n\rangle_{I_n}}{\|\psi_n\|_{I_n,2}}\left|\sqrt{\mathrm{Var}_{l_n}\widehat{S}_n}\bigg/\left(\frac{1}{\sharp I_n}\sum_{i\in I_n}\psi_n(X_i)^2\frac{R_n(i)^2}{(\sharp I_n+1)^2}\right)^{1/2}-1\right|$$
$$=o_{P_{l_n}}(1).$$

Since $\bar{G}(\cdot)$ is uniformly bounded by 1, the Lindeberg condition is easily verified for $\widehat{S}_n$. Then Lindeberg's central limit theorem yields in combination with the result from step 1, (17) and (23)

$$\mathbb{P}_{l_n}(T_{j_nk_n}>C+\sqrt{2\log(n/\sharp I_n)})$$
$$=1-\Phi\left(C+\sqrt{2\log(n/\sharp I_n)}-\sqrt{12}\sqrt{\sharp I_n}\frac{\langle\psi_n,l_n\rangle_{I_n}}{\|\psi_n\|_{I_n,2}}\int f(y)^2\,dy\right)+o(1),$$

with $\Phi$ the standard normal distribution function. It remains to be shown that

(24) $$\sqrt{12}\sqrt{\sharp I_n}\frac{\langle\psi_n,l_n\rangle_{I_n}}{\|\psi_n\|_{I_n,2}}\int f(y)^2\,dy-\sqrt{2\log(n/\sharp I_n)}\to\infty$$

as $n$ goes to infinity under the constraints $\|l_n\sqrt{h}\|_{\sup}\geq d^*\rho_n(1+\varepsilon_n)$ and (16).

Under the assumptions about the kernel $\psi$ and the design density $h$, arguments involving bounded total variation of $\psi$ and $h$ yield the approximation

(25)
$$\sqrt{12}\sqrt{\sharp I_n}\frac{\langle\psi_n,l_n\rangle_{I_n}}{\|\psi_n\|_{I_n,2}}\int f(y)^2\,dy-\sqrt{2\log(n/\sharp I_n)}$$
$$=\sqrt{12}\sqrt{n}\frac{\langle\psi_n,l_n\sqrt{h}\rangle}{\|\psi_n\|_2}\int f(y)^2\,dy-\sqrt{2\log(n/(\sharp I_n))}+o(1).$$

Let $\psi^{(n)}$ be the kernel rescaled to the interval $[M(l_n)-h_n,M(l_n)+h_n]$. Then

$$\frac{\langle\psi_n,l_n\sqrt{h}\rangle}{\|\psi_n\|_2}=\frac{\langle\psi^{(n)},l_n\sqrt{h}\rangle}{\|\psi^{(n)}\|_2}(1+O((nh_n)^{-1})),$$

using that $X_{j_n}-(M(l_n)-h_n)=O(n^{-1})$ and $X_{k_n}-(M(l_n)+h_n)=O(n^{-1})$ by assumption (D). But $\delta_n\psi^{(n)}$ by its construction as well as $l_n\sqrt{h}$ are elements of $\mathcal{H}(\beta,L)$. Then as in Dümbgen and Spokoiny (2001), a convexity argument yields the inequality

$$\delta_n^{-1} \frac{\langle \delta_n \psi^{(n)}, l_n \sqrt{h} \rangle}{\|\psi^{(n)}\|_2} \geq \frac{\delta_n^{-1} \|\delta_n \psi^{(n)}\|_2^2}{\|\psi^{(n)}\|_2} = \delta_n \sqrt{h_n} \|\gamma_\beta\|_2. \tag{26}$$

One verifies that

$$\sqrt{12}\sqrt{n}\left(\int f(y)^2\, dy\right) \delta_n \sqrt{h_n} \|\gamma_\beta\|_2 (1 + O((nh_n)^{-1})) - \sqrt{2\log(1/h_n)} + o(1)$$
$$\geq \varepsilon_n (2/(2\beta+1))^{1/2} \sqrt{\log n} + o(1) \to \infty$$

and therefore (24) follows in combination with (25) and (26).

*Step* 3. Suppose now that there exists a sequence $(l_n)$ with

$$\liminf_{n \to \infty} \mathbb{P}_{l_n}(T_{j_n k_n} > C + \sqrt{2\log(n/\sharp I_{j_n k_n})}) = c < 1,$$

where the indices $j_n, k_n$ are chosen as in step 2. This implies the existence of a subsequence [for simplicity also denoted by $(l_n)$] without any subsubsequence having the property (16); that is, we may assume $\|l_n\|_{\sup}/\rho_n \to \infty$. We will conclude the proof via contradiction as follows: For any subsequence of a sequence $(l_n)$ satisfying $\|l_n\|_{\sup}/\rho_n \to \infty$, there exists a subsubsequence which either converges to zero at a slow rate or whose maximal absolute value stays uniformly bounded away from zero. Hence we need to show that in both cases, our test attains asymptotic power 1.

Note that the squared denominator of $T_{j_n k_n}$ is bounded by $\|\psi\|_{\sup}^2$, while $\operatorname{Var}_{l_n}(\widehat{S}_n)$ is uniformly bounded. Using again the approximation of the numerator by $\widehat{S}_n$, we obtain

$$\mathbb{E}_{l_n} T_{j_n k_n} - \sqrt{2\log(n/\sharp I_n)}$$
$$\geq \|\psi\|_{\sup}^{-1} \mathbb{E}_{l_n} \widehat{S}_n - \sqrt{2\log(n/\sharp I_n)} + o(1). \tag{27}$$

If there exists a sequence $(l_n)$ with the property $\|l_n\|_{\sup}/\rho_n \to \infty$ but which converges to zero,

$$\mathbb{E}_{l_n} \widehat{S}_n = 2\sqrt{\sharp I_n} \langle \psi_n, l_n \rangle_{I_n} \left\{ \int f(y)^2\, dy \right\} + \sqrt{\sharp I_n} O(\|l_n\|_{I_n,2}^2), \tag{28}$$

as seen in step 2. But then the first term dominates in order the second one as well as the logarithmic correction which shows that the right-hand side in (27) goes to infinity.

Otherwise, assume that $(l_n)$ stays uniformly bounded away from zero. First observe that with $\bar{l}_n := 1/(\sharp I_n) \sum_{i \in I_n} l_n(X_i)$, $|l_n(X_i) - \bar{l}_n(X_i)| \leq L|X_{j_n} -$ $X_{k_n}|^\beta = O(h_n^\beta)$. Taylor expansion around $\bar{l}_n$ up to the first order provides the approximation

$$\mathbb{E}_{l_n}\widehat{S}_n = \frac{1}{\sqrt{\sharp I_n}} \sum_{i \in I_n} \psi_n(X_i) \left\{ \int (\bar{F}(y + l_n(X_i)) - \bar{F}(-y - l_n(X_i))) f(y)\, dy \right\}$$

$$= \frac{1}{\sqrt{\sharp I_n}} \sum_{i \in I_n} \psi_n(X_i) \left\{ \int (F(y) - F(-y - 2\bar{l}_n)) f(y)\, dy + O(h_n^\beta) \right\}$$

$$= \mathbb{E}_{\bar{l}_n}\widehat{S}_n + O(n^{1/2} h_n^{\beta + 1/2}).$$

If $\bar{l}_n$ is uniformly bounded away from zero, $\mathbb{E}_{\bar{l}_n}\widehat{S}_n$ is of order not smaller than $O(\sqrt{nh_n})$ which dominates in order the approximation error $|\mathbb{E}_{\bar{l}_n}\widehat{S}_n - \mathbb{E}_{l_n}\widehat{S}_n|$ as well as the logarithmic correction. $\square$

PROOF OF THEOREM 4. By virtue of the proof of Theorem 3, it remains to be shown that (i) there exists some positive constant $C = C(\beta, L, \psi)$, such that (24) goes to infinity for alternatives $l_n$ with $K\rho_n \geq \|l_n \sqrt{h}\|_{\sup} \geq C\rho_n$ for any constant $K > C$ and (ii) $\mathbb{E}_{l_n}\widehat{S}_n$ goes to infinity whenever $\|l_n\|_{\sup}/\rho_n \to \infty$. To this aim, we establish the following: If $l \in \mathcal{H}(\beta, L)$ with $\|l\|_{\sup} \leq 1$ and $x^* := \arg\max_{x \in [0,1]} |l(x)|$, then there exist some constant $c = c(\beta, L) > 0$ and a closed interval $I(l) \subset [0,1]$ such that $\lambda(I(l)) \geq c|l(x^*)|^{1/\beta}$ and

(29) $\qquad |l(x)| \geq \tfrac{1}{2}|l(x^*)| \qquad$ for every $x \in I(l)$.

Note that this is obviously correct in case $\beta \leq 1$ with $c = 1/(2L)$. For $\beta > 1$, let $\lfloor \beta \rfloor$ denote the largest integer strictly smaller than $\beta$. Let $l \in \mathcal{H}(\beta, L)$ with $\|l\|_{\sup} = D > 0$. Taylor expansion around any point $y \in [0,1]$ provides the approximation

$$l(x) = l(y) + (x-y)l'(y) + \cdots + \frac{(x-y)^{\lfloor \beta \rfloor}}{k!} l^{(\lfloor \beta \rfloor)}(y) + R_l(x, y)$$

with $|R_l(x,y)| \leq L|x-y|^\beta (\leq L)$. Thus,

(30) $\qquad \left|(x-y)l'(y) + \cdots + \frac{(x-y)^{\lfloor \beta \rfloor}}{k!} l^{(\lfloor \beta \rfloor)}(y)\right| \leq 2D + L.$

LEMMA. *There exists a universal constant $K = K_d$ such that for any polynomial $P$ of degree $d > 0$, say $P(x) = \sum_{k=0}^{d} a_k x^k$, and $\|P\|_{[0,1]} \leq D > 0$, it holds true that $\sup_{k=0,\ldots,d} |a_k| \leq K_d \cdot D$.*

The lemma results from the fact that, for the polynomial $P(x) = \sum_{k=0}^{d} a_k x^k$, $\|P\|_{(1)} = \|P\|_{[0,1]}$ and $\|P\|_{(2)} = \max_{0 \leq k \leq d} |a_k|$ are two norms in the $(d+1)$-dimensional space of polynomials of degree $d$, and these norms are equiva-lent. Its application implies together with the bound (30) that there exists a constant $K = K(\beta)$ such that $|l(x) - l(x^*)| \leq \|l'\|_{\sup}|x - x^*| \leq K(2D + L)|x - x^*|$. Then $|l(x)| \geq 1/2|l(x^*)|$ on $[x^* - D/(4KD + 2KL), x^* + D/(4KD + 2KL)] \cap [0, 1]$. If now $l_n \in \mathcal{H}(\beta, L)$ with $\|l_n\|_{\sup} = \delta_n \leq 1$, then at least $[x^* - 2^{-1}\delta_n^{1/\beta}, x^*]$ or $[x^*, x^* + 2^{-1}\delta_n^{1/\beta}]$ is fully contained in $[0, 1]$. Assume without loss of generality that $[x^*, x^* + 2^{-1}\delta_n^{1/\beta}] \subset [0, 1]$. Then $g_n$ is defined by $g_n(x) := 2^\beta \delta_n^{-1} l_n(2^{-1}\delta_n^{1/\beta}x + x^*)$ for $x \in [0, 1]$ is element of $\mathcal{H}(\beta, L)$ with $\|g_n\|_{\sup} = g_n(0) = 2^\beta$. Thus the above lemma finally implies that $|l_n(x)| \geq \delta_n/2$ on $[x^*, x^* + 1/(8K + 4K2^{-\beta}L)\delta_n^{1/\beta}]$.

The assumption about $\psi$ implies that there exists some interval $[c, d] \subset (0, 1)$ on which $\psi(x) \geq \delta$ for some strictly positive constant $\delta$. We first verify the claim (i). For any alternative $l_n$, let $\psi_n$ be the kernel rescaled onto the interval $[X_{j_n}, X_{k_n}]$, where the design points $X_{j_n} < X_{k_n}$ are those which are closest to the endpoints of $I(l_n\sqrt{h})$. Let $I_n := \{i : X_i \in I(l_n\sqrt{h})\}$. Then $\langle \psi_n, l_n\sqrt{h}\rangle_{I_n}$ is of order not smaller than $\|l_n\sqrt{h}\|_{\sup}$, which implies the existence of a universal constant $C = C(\beta, L, \psi)$ such that (24) goes to infinity for $\|l_n\sqrt{h}\|_{\sup} \geq C\rho_n$ and $\|l_n\|_{\sup}/\rho_n = O(1)$. The same consideration also shows that (28) goes to infinity whenever $\|l_n\|_{\sup}/\rho_n \to \infty$ and $\|l_n\|_{\sup} \to 0$, because $\|l_n\sqrt{h}\|_{\sup}$ dominates in order $\|l_n\|_{I_n,2}^2$ as well. To verify (ii), note that $\|l_n\sqrt{h}\|_{\sup}/(4K\|l_n\sqrt{h}\|_{\sup} + 2KL)$ stays uniformly bounded away from zero and infinity as soon as $\|l_n\|_{\sup}$ is uniformly bounded away from zero. Thus in the latter case, there always exists an interval $I(l_n\sqrt{h})$ with $\liminf_{n\to\infty} \lambda(I(l_n\sqrt{h})) > 0$ and $|l_n(X_i)\sqrt{h(X_i)}| \geq \|l_n\sqrt{h}\|_{\sup}/2$ for every $X_i \in I(l_n\sqrt{h})$. With $I_n := \{i|X_i \in I(l_n\sqrt{h})\}$

$$S_n = \frac{1}{\sqrt{\sharp I_n}} \sum_{i \in I_n} \psi_n(X_i) \operatorname{sign}(Y_i) \frac{\mathbb{E}_{l_n}(R_n(i)|\operatorname{sign}(Y_i))}{\sharp I_n + 1}$$
$$+ \frac{1}{\sqrt{\sharp I_n}} \sum_{i \in I_n} \psi_n(X_i) \operatorname{sign}(Y_i) \frac{R_n(i) - \mathbb{E}_{l_n}(R_n(i)|\operatorname{sign}(Y_i))}{\sharp I_n + 1}.$$

If $l_n(X_i)$ is uniformly bounded away from zero for every $i \in I_n$, the absolute expectation of first term is of order $O(\sqrt{n})$, while the second term is $O_p(1)$. □

**Acknowledgments.** The present work is part of my Ph.D. thesis. I sincerely thank my supervisor Lutz Dümbgen for his constant encouragement and support. Furthermore, I am grateful to the Co-Editor Jianqing Fan, an Associate Editor and two referees for many suggestions which lead to a substantial improvement of this article.

Weierstrass-Institut Für
 Angewandte Analysis und Stochastik
Mohrenstrasse 39
D-10117 Berlin
Germany
E-mail: rohde@wias-berlin.de